\documentclass[11pt,a4]{article}
\usepackage{amsmath,amssymb,amsfonts,amsthm,euscript,amsbsy}
\usepackage[dvips]{graphicx}
\usepackage{epic}

\usepackage[usenames]{color}

\newcommand{\ch}[1]{#1}

\newcommand{\wh}[1]{\widehat{#1}}
\newcommand{\wt}[1]{\widetilde{#1}}
\newtheorem{THM}{Theorem}[section]

\newtheorem{REM}[THM]{Remark}
\newtheorem{LEM}[THM]{Lemma}

\newcommand{\R}[1]{\frac{1}{#1}}

\def\Rd{\mathbb{R}^d}
\def\sp{\enspace}
\def\del{\partial}

\def\t0{\wt{t}_0(0)}

\def\hlf{\frac{1}{2}}

\def\dfour{\frac{d-4}{2}}

\def\N{\mathbb{N}}
\def\Z{\mathbb{Z}}

\renewcommand{\to}{\rightarrow}

\newcommand{\nn}{\nonumber}

\newcommand{\eqn}[1]{\begin{equation}#1\end{equation}}
\newcommand{\seqn}[1]{\begin{equation}\begin{split}#1
  \end{split}\end{equation}}
\newcommand{\eqnlab}[2]{\begin{equation}\label{#1}#2\end{equation}}
\newcommand{\seqnlab}[2]{\begin{equation}\label{#1}\begin{split}#2
  \end{split}\end{equation}}

\newcommand{\sm}[3]{\frac{#1}{#2^{#3}}}
\newcommand{\mc}[1]{\mathcal{#1}}

\numberwithin{equation}{section}

\title{Extension of the generalised inductive approach to the lace expansion:\\
Full proof}
\author{Remco van der Hofstad
\thanks{Department of Mathematics and
Computer Science, Eindhoven University of Technology, P.O.\ Box
513, 5600 MB Eindhoven, The Netherlands. E-mail:
{rhofstad@win.tue.nl, holmes@eurandom.tue.nl} }
\and Mark Holmes$\;^*$
\and Gordon Slade
\thanks{Department of Mathematics, University of
British Columbia, Vancouver, BC V6T 1Z2, Canada. E-mail:
{slade@math.ubc.ca} } }

\date{May 10, 2007}

\begin{document}
\maketitle
\begin{abstract}
This paper extends the result of \cite{HS02} in order to use the inductive approach to prove Gaussian asymptotic behaviour for models with critical dimension other than $4$.  The results are applied in \cite{H07} to study sufficiently spread-out lattice trees in dimensions $d>8$ and may also be applicable to percolation in dimensions $d>6$.
\end{abstract}

\section{Introduction}
This paper consists of large parts of the material in \cite{HS02}, reproduced verbatim, but with the introduction of parameters $\theta(d)>2$ and $p^*\ge 1$ as described in \cite{HHS07}.  The case $\theta=\frac{d}{2}$, and $p^*=1$ is that dealt with in \cite{HS02}.  The main conclusion is one of Gaussian asymptotic behaviour for models with critical dimension other than $4$, satisfying certain properties.  We do not include the proof of the local central
limit theorem \cite[Theorem~1.3]{HS02}, which does require $\theta = \frac d2$.
The result of this paper is applied in \cite{H07} to lattice trees with $d>8$, $\theta=\dfour$ and $p^*=2$.
We also expect the result to be applicable to other models where the analysis uses the lace expansion above a critical dimension $d_c\ge 4$.  In such cases the lace expansion for $d>d_c$ suggests setting $\theta=\frac{d-(d_c-4)}{2}$.
In particular the above statement for percolation in dimensions $d>d_c=6$ would give $\theta=\frac{d-2}{2}$.

This paper simply provides the details of the proof described in \cite{HHS07}, and we refer the reader to \cite{HHS07} and \cite{HS02} for a more thorough introduction to the inductive approach to the lace expansion.  In Section \ref{sec:assthm} we state the form of the recursion relation, and the assumptions \ch{S, D, E$_{\theta}$, and G$_{\theta}$}
on the quantities appearing in the recursion equation.  We also state the ``$\theta$-theorem'' to be proved.  In
Section~\ref{sec-ih}, we introduce the induction hypotheses on
$f_n$ that will be used to prove the $\theta$-theorem, and derive some
consequences of the induction hypotheses. The induction is advanced in
Section~\ref{sec-adv}. In Section~\ref{sec-pf}, the $\theta$-theorem
stated in Section~\ref{sec:assthm} are proved.
\section{Assumptions on the Recursion Relation}
\label{sec:assthm}
When applied to self-avoiding walks, oriented percolation and lattice trees, the
lace expansion gives rise to a convolution recursion relation of
the form
        \eqnlab{fkrec}{
        f_{n+1}(k;z) =
        \sum_{m=1}^{n+1} g_m(k;z)f_{n+1-m}(k;z) + e_{n+1}(k;z)
        \quad \quad (n \geq 0),}
with $f_0(k;z) = 1$.
Here, $k \in [-\pi,\pi]^d$ is a parameter dual to a
spatial lattice variable $x \in \Z^d$, and $z$ is a positive parameter.
The functions $g_m$ and $e_m$ are to be regarded as given, and the goal
is to understand the behaviour of the solution $f_n(k;z)$ of (\ref{fkrec}).

\subsection{Assumptions S,D,E$_\theta$,G$_\theta$}
The first assumption,
Assumption S, requires that the functions appearing in the recursion equation (\ref{fkrec}) respect the lattice symmetries of reflection and rotation,
and that $f_n$ remains bounded in a weak sense.  We have strengthened this assumption from that appearing in \cite{HS02}, as one requires smoothness of $f_n$ and $g_n$ which holds in all of the applications.

\smallskip \noindent
{\bf Assumption S}. For every $n\in\N$ and $z>0$, the mapping
$k\mapsto f_n(k;z)$ is symmetric under replacement of any component
$k_i$ of $k$ by
$-k_i$, and under permutations of the components of $k$.  The same
holds for $e_n(\cdot;z)$ and $g_n(\cdot; z)$.  In addition, for
each $n$, $|f_n(k;z)|$ is bounded uniformly in $k \in
[-\pi,\pi]^d$ and $z$ in a neighbourhood of $1$ (which may depend on $n$).  We also assume that $f_n$ and $g_n$ have continuous second derivatives in a neighbourhood of $0$ for every $n$.  It is an immediate consequence of Assumption S that the mixed partials of $f_n$ and $g_n$ at $k=0$ are equal to zero.

\smallskip
The next assumption, Assumption~D,
incorporates a ``spread-out'' aspect to the recursion equation.  It introduces a function $D$ which defines the underlying random walk model, about which Equation (\ref{fkrec}) is a perturbation.
The assumption involves a non-negative parameter $L$, which will
be taken to be large, and which serves to spread out the steps of the
random walk over a large set.
We write $D=D_L$ in the statement of Assumption~D to emphasise this
dependence, but the subscript will not be retained elsewhere.
An example of a family of $D$'s obeying the assumption is taking $D(\cdot)$ uniform on a box side length $2L$, centred at the origin.  In particular Assumption~D implies that $D$ has a finite second moment and we define
\eqnlab{sigdef}{\sigma^2 \equiv  - \nabla^2 \hat{D}(0)=-\left[\sum_j \frac{\del^2 }{\del k_j^2}\sum_x e^{ik\cdot x}D(x)\right]_{k=0}=-\left[\sum_j \sum_x (ix_j)^2 e^{ik\cdot x}D(x)\right]_{k=0}=\sum_x |x|^2 D(x).}

The assumptions involve a parameter $d$, which corresponds to the spatial dimension in our applications, and a parameter $\theta> 2$ which will be model dependent.

Let
\eqnlab{adef}{a(k) = 1 - \hat{D}(k).}
\smallskip \noindent
{\bf Assumption D}.
We assume that $D(x)\ge 0$ and
\eqn{f_1(k;z) = z \hat{D}_L(k), \quad e_1(k;z)=0.}
In particular, this implies that
$g_1(k;z)=z\hat{D}_L(k)$.  As part of Assumption~D, we also assume:\\
(i)
$D_L$ is normalised so that
$\hat{D}_L(0) =1$, and has $2+2\epsilon$ moments for some
$\epsilon \in (0,\theta-2)$, i.e.,
\eqnlab{momentD}{\sum_{x\in \Z^d} |x|^{2+2\epsilon} D_L(x) <\infty.}
(ii)
There is a constant $C$ such that, for all $L \geq 1$,
\eqnlab{beta, sigmadef}{\|D_L\|_\infty \leq CL^{-d},
\qquad \sigma^2 = \sigma^2_L\leq CL^2,}\\
(iii)
There exist
constants $\eta,c_1,c_2 >0$ such that
\eqnlab{Dbound1}{c_1 L^2 k^2 \leq a_L(k) \leq c_2 L^2 k^2 \quad (\|k\|_\infty \leq L^{-1}),}
\eqnlab{Dbound2}{a_L(k) > \eta  \quad (\|k\|_\infty \geq L^{-1}),}
\eqnlab{Dbound3}{a_L(k) < 2-\eta \quad (k \in [-\pi,\pi]^d).}


Assumptions~E and~G of \cite{HS02} are now adapted to general $\theta>2$ as follows.
The relevant bounds on $f_m$, which {\em a priori}\/ may or may
not be satisfied, are that for some $p^*\ge 1$, some nonempty $B\subset[1,p^*]$ and
\eqn{\beta=\beta(p^*)=L^{-\frac{d}{p^*}}}
 we have for every $p \in B$,
\eqnlab{fbdsp}{
\|\hat{D}^2 f_m(\cdot;z)\|_p\leq \frac{K}{L^{\frac{d}{p}} m^{\frac{d}{2p} \wedge \theta}},
    \quad  | f_m(0;z)|\leq K, \quad
    |\nabla^2 f_m(0;z)|\leq K \sigma^2 m,}
for some positive constant $K$.  The full generality in which this has been presented is not required for our application to lattice trees where we have $p^*=2$ and $B=\{2\}$.  This is because we require only the $p=2$ case in (\ref{fbdsp}) to estimate the diagrams arising from the lace expansion for lattice trees and verify the assumptions \ch{E$_{\theta}$, G$_{\theta}$} which follow.  In other applications it may be that a larger collection of
\ch{$\|\cdot\|_p$} norms are required to verify the assumptions and the set $B$ is allowing for this possibility.  The parameter $p^*$ serves to make this set bounded so that $\beta(p^*)$ is small for large $L$.

The bounds in (\ref{fbdsp}) are
identical to the ones in \cite{HS02}, except for the first bound,
which only appears for $p=1$ and $\theta=\frac{d}{2}$.

\smallskip \noindent
{\bf Assumption E$_{\theta}$}. There is an $L_0$, an interval $I \subset
[1-\alpha,1+\alpha]$ with $\alpha \in (0,1)$, and a function $K
\mapsto C_e(K)$, such that if (\ref{fbdsp}) holds for some $K>1$, $L \geq L_0$, $z \in I$ and for all $1 \leq m \leq n$,
then for that $L$ and $z$, and for all $k \in [-\pi,\pi]^d$ and $2
\leq m\leq n+1$, the following bounds hold:
\eqn{|e_m(k;z)|\leq C_e(K) \beta m^{-\theta},
    \quad |e_m(k;z)-e_m(0;z)|\leq
    C_e(K) a(k) \beta m^{-\theta+1}.}

\smallskip \noindent
{\bf Assumption G$_{\theta}$}. There is an $L_0$, an interval $I \subset
[1-\alpha,1+\alpha]$ with $\alpha \in (0,1)$, and a function $K
\mapsto C_g(K)$, such that if (\ref{fbdsp}) holds for some $K>1$, $L \geq L_0$, $z \in I$ and for all $1 \leq m \leq n$, then for
that $L$ and $z$, and for all $k \in [-\pi,\pi]^d$ and $2 \leq
m\leq n+1$, the following bounds hold:
\eqn{
    |g_m(k;z)|\leq C_g(K)\beta m^{-\theta},
    \quad |\nabla^2 g_m(0;z)|\leq
    C_g(K) \sigma^2 \beta m^{-\theta+1},}
\eqn{|\partial_z g_m(0;z)|\leq
    C_g(K)\beta m^{-\theta+1},}
\eqn{|g_m(k;z)-g_m(0;z)- a(k) \sigma^{-2}
    \nabla^2 g_m(0;z)| \leq C_g(K)\beta
    a(k)^{1+\epsilon'}m^{-\theta+(1+\epsilon')},}
with the last bound valid for any $\epsilon' \in [0,\epsilon]$.
\begin{THM}
\label{thm-1p} Let $d>d_c$ and $\theta(d)>2$,
and assume that Assumptions $S$, $D$,
$E_{\theta}$ and $G_{\theta}$ all hold. There exist positive $L_0 = L_0(d,\epsilon)$,
$z_c=z_c(d,L)$, $A=A(d,L)$, and $v = v(d,L)$, such that for $L
\geq L_0$, the following statements hold.
\\
(a)  Fix $\gamma  \in (0,1 \wedge \epsilon)$
and $\delta \in (0, (1\wedge \epsilon) -
\gamma)$.  Then
\eqn{f_n\Big(\frac{k}{\sqrt{v \sigma^2 n}};z_c\Big) =A e^{-\frac{k^2}{2d}}
        [1 + {\cal O}(k^2n^{-\delta})+{\cal O}(n^{-\theta+2})],}
with the error estimate uniform in $\{k \in \Rd: a(k/\sqrt{v\sigma^2 n}) \leq \gamma n^{-1} \log n \}$.
\\
(b)
\eqn{-\frac{\nabla^2 f_n(0;z_c)}{f_n(0;z_c)}=v \sigma^2 n [1+{\cal O}(\beta n^{-\delta})].}
(c)  For all $p\ge1$,
\eqn{\|\hat{D}^2 f_n(\cdot;z_c)\|_p\leq \frac{C}{L^{\frac{d}{p}}n^{\frac{d}{2p}\wedge \theta}}.}
(d) The constants $z_c$, $A$ and $v$ obey
\seqnlab{eq:pthmd}{1 &=  \sum_{m = 1}^\infty g_m(0;z_c),\\
    A &= \frac{1+\sum_{m=1}^\infty e_m(0;z_c)}{\sum_{m=1}^\infty m g_m(0;z_c)},\\
    v &=-\frac{\sum_{m=1}^\infty \nabla^2 g_m(0;z_c)}{\sigma^2\sum_{m = 1}^\infty m g_m(0;z_c)}.}
\end{THM}
It follows immediately from Theorem~\ref{thm-1p}(d) and the bounds of Assumptions~\ch{E$_{\theta}$ and G$_{\theta}$}
that
\eqn{z_c=1+ \mc{O}(\beta), \quad A=1+ \mc{O}(\beta), \quad v = 1+ \mc{O}(\beta).}

With modest additional assumptions, the critical point $z_c$ can be characterised in terms of the {\em susceptibility}\/
\eqnlab{sus1}{\chi(z) = \sum_{n=0}^\infty f_n(0;z).}

\begin{THM}
\label{thm-zc}
Let $d>d_c$, $\theta(d)>2$, $p^*\ge 1$ and assume that Assumptions $S$, $D$, $E_{\theta}$ and $G_{\theta}$ all hold.  Let $L$ be sufficiently large.  Suppose there is a $z_c'>0$ such that the susceptibility (\ref{sus1}) is
absolutely convergent for $z \in (0,z_c')$, with $\lim_{z \uparrow z_c'}\chi(z) = \infty$ (if $\chi(z)$ is a power series in $z$ then $z_c'$ is the radius of convergence of $\chi(z)$).  Suppose also that the bounds of (\ref{fbdsp}) for $z=z_c$ and all $m \geq 1$ imply the bounds of Assumptions~$E_{\theta}$ and $G_{\theta}$ for all $m \geq 2$, uniformly in $z \in [0,z_c]$.  Then $z_c=z_c'$.
\end{THM}

\section{Induction hypotheses} \label{sec-ih}

We will analyse the recursion relation (\ref{fkrec}) using induction on $n$, as done in \cite{HS02}.
In this section, we introduce the induction hypotheses, verify that
they hold for $n=1$, discuss their motivation, and derive some of
their consequences.

\subsection{Statement of induction hypotheses (H1--H4)}
    \label{sec-ihstate}
The induction hypotheses involve a sequence $v_n$, which is
defined as follows.  We set $v_0=b_0=1$, and for $n \geq
1$ we define
    \eqn{
        b_n = -\frac{1}{\sigma^2}\sum_{m=1}^{n}
    \nabla^2 g_m(0;z)
    ,\quad
    c_n  =  \sum_{m=1}^{n} (m-1) g_m(0;z),\quad
    \label{Delta_n}
    v_n  = \frac{b_n}{1+c_n}.
    }
The $z$--dependence of $b_n$, $c_n$, $v_n$ will usually be left implicit
in the notation.  We will often simplify the notation
by dropping $z$ also from $e_n$, $f_n$ and $g_n$, and write, e.g.,
$f_n(k)=f_n(k;z)$.

\begin{REM}
\label{rem:b1}
Note that the above definition and assumption~D gives
\eqnlab{eq:b1}{b_1=-\frac{1}{\sigma^2}\nabla^2g_1(0;z)=-\frac{1}{\sigma^2}\nabla^2 z\wh{D}(0)=-\frac{z}{\sigma^2}.\left(-\sigma^2\right)=z.}
Obviously we also have $c_1=0$ so that $v_1=z$.
\end{REM}


The induction hypotheses also involve several constants. Let $d>d_c$, $\theta>2$,
and recall that $\epsilon$ was specified in (\ref{momentD}).  We fix
$\gamma, \delta>0$ and $\lambda>2$ according to
\seqnlab{agddef}{0<\gamma&<1 \wedge \epsilon\\
0<\delta&<(1 \wedge\epsilon)-\gamma\\
\theta-\gamma<\lambda&<\theta.}



\noindent We also introduce constants $K_1, \ldots , K_5$, which are independent of $\beta$.
We define
\eqnlab{K4'def}{
        K_4' = \max \{C_e(cK_4), C_g(cK_4), K_4\},
}
where $c$ is a constant determined in Lemma~\ref{lem-pibds} below.
To advance the induction, we will need to assume that
    \eqnlab{Kcond}{
    K_3 \gg K_1 > K_4' \geq K_4 \gg 1, \quad
    K_2 \geq K_1, 3K_4' , \quad
    K_5 \gg K_4.
    }
Here $a \gg b$ denotes the statement that $a/b$ is sufficiently large.
The amount by which, for instance, $K_3$ must exceed $K_1$ is
independent of $\beta$ (but may depend on $p^*$) and will be determined during the course of
the advancement of the induction in Section~\ref{sec-adv}.

Let $z_0=z_1=1$, and define $z_n$ recursively by
    \eqn{
    \label{z_n}
    z_{n+1} = 1-\sum_{m=2}^{n+1}g_m(0;z_n),
    \qquad n \geq 1.
    }
For $n \geq 1$, we define intervals
    \eqnlab{Indef}{
    I_n = [z_n - K_1\beta n^{-\theta+1}, z_n + K_1\beta n^{-\theta+1}].
    }
In particular this gives $I_1=[1-K_1\beta, 1+K_1\beta]$.

Recall the definition $a(k)=1-\hat{D}(k)$ from (\ref{adef}). Our
induction hypotheses are that the following four statements hold
for all $z \in I_n$ and all $1\leq j\leq n$.

\begin{description}
\item[(H1)]
$|z_j - z_{j-1}| \leq K_1 \beta j^{-\theta}$.
\item[(H2)]
$|v_j - v_{j-1}| \leq K_2 \beta j^{-\theta+1}$.
\item[(H3)]
For $k$ such that $a(k) \leq \gamma j^{-1}\log j$, $f_j(k;z)$ can
be written in the form
\[
    f_j(k;z) = \prod_{i=1}^j\left[
    1 -v_i a(k)+ r_i(k)  \right],
\]
with $r_i(k) = r_i(k;z)$ obeying
\[
    |r_i(0)|\leq K_3 \beta i^{-\theta+1},\quad
    |r_i(k)-r_i(0)| \leq K_3 \beta a(k) i^{-\delta}.
\]
\item[(H4)]
For $k$ such that $a(k) > \gamma j^{-1}\log j$, $f_j(k;z)$ obeys
the bounds
\[
    |f_j(k;z)| \leq K_4 a(k)^{-\lambda}j^{-\theta}, \quad
    |f_j(k;z) -f_{j-1}(k;z)| \leq K_5 a(k)^{-\lambda+1} j^{-\theta}.
\]
\end{description}

\medskip \noindent
Note that, for $k=0$, (H3) reduces to $f_j(0) = \prod_{i=1}^j [1 + r_i(0) ]$.

\subsection{Initialisation of the induction}

We now verify that the induction hypotheses hold when $n=1$.  This remains unchanged from the $p=1$ case.  Fix $z \in I_1$.

\begin{description}
\item[(H1)]
We simply have $z_1-z_0=1-1=0$.
\item[(H2)]
From Remark \ref{rem:b1} we simply have $|v_1-v_0| = |z-1|$, so that (H2) is satisfied provided
$K_2 \ge K_1$.
\item[(H3)]
We are restricted to $a(k)=0$.  By (\ref{Dbound1}),
this means $k=0$. By
Assumption~D, $f_1(0;z) = z$, so that $r_1(0) = z-1=z-z_1$.  Thus (H3) holds provided we take $K_3 \ge K_1$.
\item[(H4)]
We note that $|f_1(k;z)| \leq z \leq 2$ for $\beta$
sufficiently small (i.e. so that $\beta K_1 \le 1$), $|f_1(k;z)-f_0(k;z)| \leq 3$, and $a(k) \leq 2$.  The bounds of (H4) therefore hold provided we take $K_4 \ge 2^{\lambda+1}$ and $K_5 \ge 3 \cdot 2^{\lambda-1}$.
\end{description}

\subsection{Discussion of induction hypotheses}
\label{sec-mot}

\noindent {\bf (H1) and the critical point.}\ The critical point
can be formally identified as follows. We set $k=0$ in (\ref{fkrec}), then
sum over $n$, and solve for the susceptibility
\eqnlab{sus}{\chi(z)=\sum_{n=0}^{\infty}f_n(0;z).}
The result is
    \eqnlab{chiz}{
    \chi(z)
    = \frac{1+\sum_{m=2}^{\infty} e_m(0;z)}{1-\sum_{m=1}^\infty g_m(0;z)}.
    }
The critical point should correspond to the smallest zero of the
denominator and hence should obey the equation
    \eqnlab{1pi}{
        1 - \sum\limits_{m=1}^{\infty} g_m(0; z_c)
    = 1 - z_c - \sum\limits_{m=2}^{\infty} g_m(0; z_c) =0.
    }
However, we do not know {\em a priori}\/ that the series in (\ref{chiz}) or
(\ref{1pi}) converge. We therefore approximate (\ref{1pi}) with
the recursion (\ref{z_n}), which bypasses the convergence issue by
discarding the $g_m(0)$ for $m>n+1$ that cannot be handled at the
$n^{\rm th}$ stage of the induction argument. The sequence $z_n$
will ultimately converge to $z_c$.

In dealing with the sequence $z_n$, it is convenient to formulate
the induction hypotheses for a small interval $I_n$ approximating
$z_c$.  As we will see in Section \ref{sec-prel}, (H1) guarantees
that the intervals $I_j$ are decreasing: $I_1 \supset I_2 \supset
\cdots \supset I_n$. Because the length of these intervals is
shrinking to zero, their intersection $\cap_{j=1}^\infty I_j$ is a
single point, namely $z_c$. Hypothesis (H1) drives the convergence
of $z_n$ to $z_c$ and gives some control on the rate. The rate is
determined from (\ref{z_n}) and the ansatz that the difference
$z_j-z_{j-1}$ is approximately $-g_{j+1}(0,z_c)$, with
$|g_j(k;z_c)|= \mc{O}( \beta j^{-\theta})$ as in  Assumption~G.

\subsection{Consequences of induction hypotheses}
\label{sec-prel}
In this section we derive important
consequences of the induction hypotheses.
The key result is that the induction hypotheses imply (\ref{fbdsp})
for all $1 \leq m \leq n$, from which the bounds of Assumptions~\ch{E$_{\theta}$
and G$_{\theta}$} then follow, for $2 \leq m \leq n+1$.

Here, and throughout the
rest of this paper:
    \begin{itemize}
    \item $C$ denotes a strictly positive constant that may depend
    on $d,\gamma,\delta,\lambda$, but {\it not}\/ on the $K_i$,
    {\it not}\/ on $k$, {\it not}\/ on $n$, and {\it not}\/ on
    $\beta$ (provided $\beta$ is sufficiently small, possibly
    depending on the $K_i$). The value of $C$ may change
    from line to line.
    \item We frequently assume $\beta \ll 1$ without explicit comment.
    \end{itemize}
The first lemma shows that the intervals $I_j$ are nested,
assuming (H1).

\begin{LEM}
\label{lem-In} Assume (H1) for $1 \leq j \leq n$. Then $I_1
\supset I_2 \supset \cdots \supset I_{n}$.
\end{LEM}

\proof
Suppose $z \in I_j$, with $2 \leq j \leq n$.  Then by
(H1) and (\ref{Indef}),
    \eqn{
    |z-z_{j-1}| \leq |z-z_j| + |z_j - z_{j-1}|
    \leq \frac{K_1 \beta}{j^{\theta-1}}
    + \frac{K_1 \beta}{j^\theta}
    \leq \frac{K_1 \beta}{(j-1)^{\theta-1}},
    }
and hence $z \in I_{j-1}$.  Note that here we have used the fact that
\seqn{\R{j^a}+\R{j^b}&
\le \R{(j-1)^{a}}\iff 1+ \R{j^{b-a}}\le \left(\frac{j}{j-1}\right)^a
}
which holds if $a\ge 1$ and $b-a\ge1$ since then
\eqn{1+ \R{j^{b-a}}\le 1+ \R{j}\le 1+ \R{j-1}\le \left(1+\frac{1}{j-1}\right)^a.}
\qed

By Lemma~\ref{lem-In}, if $z \in I_j$ for $1 \leq j \leq n$, then
$z \in I_1$ and hence, by (\ref{Indef}),
    \eqnlab{znear1}{
    |z -1| \leq K_1 \beta .
    }
It also follows from (H2)
that, for $z \in I_n$ and $1 \leq j \leq n$,
    \eqnlab{vnear1}{
    |v_j -1| \leq CK_2 \beta.
    }
Define
\eqnlab{sdef}{
        s_i(k) = [1+r_i(0)]^{-1} [v_i a(k) r_i(0) + (r_i(k) - r_i(0))].}
We claim that the induction hypothesis (H3) has the useful alternate form
\eqnlab{fs}{
        f_j(k) = f_j(0)\prod_{i=1}^j
        \left[
        1 - v_i a(k) + s_i(k)
        \right].
}
Firstly $f_j(0)=\prod_{i=1}^j [1+r_i(0)]$.  Therefore the RHS of (\ref{fs}) is
\eqn{\prod_{i=1}^j \left(1 - v_i a(k)\right)[1+r_i(0)]+v_i a(k) r_i(0) + (r_i(k) - r_i(0))}
which after cancelling terms gives the result.
Note that (\ref{fs}) shows that the $s_i(k)$ are symmetric with continuous second derivative in a neighbourhood of $0$ (since each $f_i(k)$ and $a(k)$ have these properties).  To see this note that $f_1(k)$ and $a(k)$ symmetric implies that $s_1(k)$ is symmetric.  Next, $f_2(k),a(k)$, and $s_1(k)$ symmetric implies that $s_2(k)$ symmetric etc.

We further claim that
\eqnlab{sbd}{
        |s_i(k)| \leq K_3 (2+C(K_2+K_3)\beta) \beta a(k) i^{-\delta}.
}
This is different to that appearing in \cite[(2.19)]{HS02} in that the constant is now 2 rather than 1.  This is a correction to \cite[(2.19)]{HS02} but it does not affect the analysis.  To verify (\ref{sbd}) we use the fact that $\R{1-x}\le 1+2x$ for $x \le \hlf$ to write for small enough $\beta$,
\seqn{|s_i(k)|&\le \left[1+2K_3\beta\right]\left[(1+|v_i-1|)a(k)r_i(0)+|r_i(k)-r_i(0)|\right]\\
&\le \left[1+2K_3\beta\right]\left[(1+CK_2\beta)a(k)\sm{K_3\beta}{i}{\theta-1}+\sm{K_3\beta a(k)}{i}{\delta}\right]\\
&\le \sm{K_3\beta a(k)}{i}{\delta}[1+2K_3\beta][2+CK_2\beta]\le \sm{K_3\beta a(k)}{i}{\delta}[2+C(K_2+K_3)\beta].}
Where we have used the bounds of (H3) as well as the fact that $\theta-1>\delta$.  The next lemma provides an important upper bound on $f_j(k;z)$, for $k$ small depending on $j$, as in (H3).

\begin{LEM}
\label{lem-cA} Let $z\in I_n$ and assume (H2--H3) for $1 \leq
j \leq n$. Then for $k$ with $a(k) \leq \gamma j^{-1}\log j$,
\eqn{
        |f_j(k;z)| \leq e^{CK_3\beta} e^{-(1-C(K_2+K_3)\beta)ja(k)}.
}
\end{LEM}

\proof
We use H3, and conclude from the bound on $r_i(0)$ of (H3)
that
\[|f_j(0)| = \prod_{i=1}^j |1+r_i(0)| \le \prod_{i=1}^j \left|1+\sm{K_3\beta}{i}{\theta-1}\right|\leq e^{CK_3\beta},\]
using $1+x \leq e^x$ for each
factor.
Then we use (\ref{vnear1}), (\ref{fs}) and
(\ref{sbd}) to obtain
\eqn{
        \prod_{i=1}^j  \left| 1 - v_i a(k) + s_i(k) \right|
        \leq \prod_{i=1}^j  \left| 1 - (1-CK_2\beta) a(k) +
        CK_3 \beta a(k) i^{-\delta} \right|.
}
The desired bound then follows, again using $1+x \leq e^x$ for each
factor on the right side, and by (\ref{fs}).
\qed

The middle bound of (\ref{fbdsp}) follows, for $ 1 \leq m \leq n$
and $z \in I_m$, directly
from Lemma~\ref{lem-cA}.  We next prove
two lemmas which provide the other two bounds of (\ref{fbdsp}).
This will supply the hypothesis (\ref{fbdsp}) for Assumptions~\ch{E$_{\theta}$ and G$_{\theta}$},
and therefore plays a crucial role in advancing the induction.

\begin{LEM}
\label{lem-Lpnorm}
Let $z \in I_n$ and assume (H2), (H3) and (H4).
Then for all $1 \leq j \leq n$, and $p \ge 1$,
\eqn{
        \| \hat{D}^2 f_j(\cdot ;z)\|_p \leq \frac{C(1+K_4)}{L^{\frac{d}{p}}j^{\frac{d}{2p}\wedge \theta}},}
where the constant $C$ may depend on $p,d$.
\end{LEM}

\proof We show that
\eqn{
        \| \hat{D}^2 f_j(\cdot ;z)\|_p^p \leq \frac{C(1+K_4)^p}{L^{d} j^{\frac{d}{2}\wedge \theta p}}.}
For $j=1$ the result holds since $|f_1(k)|=|z\wh{D}(k)|\le z\le 2$ and by using (\ref{beta,  sigmadef}) and the fact that $p\ge 1$.  We may therefore assume that $j\ge 2$ where needed in what follows, so that in particular $\log j \ge \log 2$.

Fix $z \in I_n$ and $1 \leq j \leq n$, and define
\seqn{
    R_{1} & = \{k \in
    [-\pi,\pi]^{d}: a(k) \leq \gamma j^{-1}\log j , \;
    \|k\|_\infty \leq L^{-1} \} ,\nonumber \\
    R_{2} & = \{k \in
    [-\pi,\pi]^{d}: a(k) \leq \gamma j^{-1}\log j , \;
    \|k\|_\infty > L^{-1} \} ,\nonumber \\
    \label{Rjdef}
    R_{3} & = \{k \in
    [-\pi,\pi]^{d}: a(k) > \gamma j^{-1}\log j, \;
    \|k\|_\infty \leq L^{-1} \},\nn\\
    R_4 & = \{k \in
    [-\pi,\pi]^{d}: a(k) > \gamma j^{-1}\log j, \;
    \|k\|_\infty > L^{-1} \}.}
The set $R_2$ is empty if $j$ is sufficiently large.  Then
    \eqn{
    \|\hat{D}^2f_{j}\|_p^p = \sum_{i=1}^4\int_{R_i}\left(\hat{D}(k)^2|f_{j}(k)|\right)^p
    \frac{d^dk}{(2\pi)^d}.
    }
We will treat each of the four terms on the right side separately.

On $R_1$, we use (\ref{Dbound1}) in conjunction with
Lemma~\ref{lem-cA} and the fact that \ch{$\hat{D}(k)^2
\leq 1$}, to obtain for all $p>0$,
    \eqn{
    \int_{R_1} \left(\hat{D}(k)^2\right)^p|f_{j}(k)|^p \frac{d^dk}{(2\pi)^d}
    \leq \int_{R_1} Ce^{-cpj(Lk)^2} \frac{d^dk}{(2\pi)^d}\le \prod_{i=1}^d \int_{-\R{L}}^{\R{L}}Ce^{-cpj(Lk_i)^2}dk_i
    \leq \frac{C}{L^d (pj)^{d/2}}\le \frac{C}{L^d j^{d/2}}.
}
Here we have used the substitution $k'_i=Lk_i\sqrt{pj}$.
On $R_2$, we use Lemma~\ref{lem-cA} and (\ref{Dbound2}) to
conclude that for all $p>0$, there is an $\alpha(p) > 1$ such that
    \eqn{
    \int_{R_2} \left(\hat{D}(k)^2|f_{j}(k)|\right)^p \frac{d^dk}{(2\pi)^d}
    \leq C\int_{R_2} \alpha^{-j} \frac{d^dk}{(2\pi)^d}
    = C\alpha^{-j} |R_2|,
}
where $|R_2|$ denotes the volume of $R_2$.  This volume is maximal
when $j=3$, so that

    \eqn{
    |R_2| \leq |\{ k : a(k)\le \textstyle\frac{\gamma \log 3}{3} \}|
    \leq |\{k: \hat{D}(k) \geq 1-\textstyle\frac{\gamma \log 3}{3} \}|
    \leq (\textstyle\frac{1}{1-\frac{\gamma \log 3}{3}})^2 \|\hat{D}^2\|_1
    \leq (\textstyle\frac{1}{1-\frac{\gamma \log 3}{3}})^2 CL^{-d},
}
using (\ref{beta, sigmadef}) in the last step.
Therefore $\alpha^{-j}|R_2| \leq C L^{-d} j^{-d/2}$ since $\alpha^{-j}j^{\frac{d}{2}}\le C(\alpha,d)$ for every $j$ (using \ch{L'H\^opital's} rule for example with $\alpha^j=e^{j\log \alpha}$), and
    \eqn{
    \int_{R_2} \left(\hat{D}(k)^2|f_{j}(k)|\right)^p \frac{d^dk}{(2\pi)^d}
    \leq CL^{-d} j^{-d/2}.
    }

On $R_3$ and $R_4$, we use (H4).  As a result, the contribution
from these two regions is bounded above by
    \eqn{
    \left(\frac{K_4 }{j^\theta}\right)^p \sum_{i=3}^4\int_{R_i}
        \frac{\hat{D}(k)^{2p}}{a(k)^{\lambda p}} \frac{d^dk}{(2\pi)^d}.}
On $R_3$, we use $\hat{D}(k)^2\leq 1$ and (\ref{Dbound1}).  Define $R_3^C=\{k:\|k\|_\infty<L^{-1}, \sp |k|^2>Cj^{-1}\log j\}$ to obtain the upper bound
    \seqn{
    \label{largekext}
    \frac{CK_4^p}{j^{\theta p} L^{2\lambda p}}
    \int_{R_3} \frac{1}{|k|^{2\lambda p}} d^dk
&\le \frac{CK_4^p}{j^{\theta p} L^{2\lambda p}}\int_{R_3^C} \frac{1}{|k|^{2\lambda p}} d^dk\\
&=\frac{CK_4^p}{j^{\theta p} L^{2\lambda p}}\int_{\sqrt{\frac{C\log j}{L^2j}}}^{\frac{d}{L}}r^{d-1-2\lambda p}dr.}
Since $\log 1 =0$, this integral will not be finite if both $j=1$ and $p \ge \frac{d}{2\lambda}$, but recall that we can restrict our attention to $j\ge 2$.  Thus we have an upper bound of
\eqn{\frac{CK_4^p}{j^{\theta p} L^{2\lambda p}}\cdot \begin{cases}
\int_0^{\frac{d}{L}}r^{d-1-2\lambda p}dr,&d>2\lambda p\\
\int_{\sqrt{\frac{C\log j}{L^2j}}}^{\frac{d}{L}}\R{r}dr,& d=2\lambda p\\
\int_{\sqrt{\frac{C\log j}{L^2j}}}^{\infty}r^{d-1-2\lambda p}dr,& d<2\lambda p
\end{cases}
\le \frac{CK_4^p}{j^{\theta p} L^{2\lambda p}}\cdot\begin{cases}
\left(\frac{d}{L}\right)^{d-2\lambda p}&, d>2\lambda p\\
\log \left(\frac{d\sqrt{L^2j}}{CL\sqrt{\log j}}\right)=\hlf\log\left(\frac{C'j}{\log j}\right)&, d=2\lambda p\\
\left(\frac{C'L^2j}{\log j}\right)^{\frac{2\lambda p-d}{2}}&, d<2\lambda p.
\end{cases}}
Now use the fact that $\lambda<\theta$ to see that each term on the right is bounded by $\frac{CK_4^p}{j^{\frac{d}{2}}L^d}$.

On $R_4$, we
use (\ref{beta, sigmadef}) and (\ref{Dbound2}) to obtain the bound
    \eqn{
    \frac{CK_4^p}{j^{\theta p}}
    \int_{[-\pi,\pi]^d}  \hat{D}(k)^{2p}\frac{d^dk}{(2\pi)^d}
    \leq \frac{CK_4^p}{j^{\theta p}} \int_{[-\pi,\pi]^d}  \hat{D}(k)^{2}\frac{d^dk}{(2\pi)^d}\le \frac{CK_4}{j^{\theta p}L^d},}
where we have used the fact that $p\ge1$ and \ch{$\hat{D}(k)^2
\leq 1$}.  Since $K_4^p \le (1+K_4)^p$, this completes the proof.
\qed

\begin{LEM}
\label{lem-fder}
Let $z \in I_n$ and assume (H2) and (H3).  Then, for
$1 \leq j \leq n$,
\eqn{
        | \nabla^2 f_j(0 ;z) | \leq (1+C(K_2 + K_3) \beta ) \sigma^2 j.
}
\end{LEM}

\proof
Fix $z\in I_n$ and $j$ with $1 \leq j \leq n$.
Using the product rule multiple times and the symmetry of all of the quantities in (\ref{fs}) to get cross terms equal to $0$,
    \eqnlab{1.2b2}{
    \nabla^2 f_{j}(0)
    = f_j(0)\sum_{i=1}^{j} \bigl[-\sigma^2v_i+\nabla^2 s_i(0)\bigr].
    }
By (\ref{vnear1}), $|v_i -1| \leq CK_2\beta$.  For the second term on the
right side, we let
$e_1,\ldots,e_d$ denote the standard basis vectors in
$\mathbb{R}^d$.  Since $s_i(k)$ has continuous second derivative in a neighbourhood of $0$, we use the extended mean value theorem $s(t)=s(0)+ts'(0)+\hlf t^2 s''(t^*)$ for some $t^*\in (0,t)$, together with (\ref{sbd}) to see that for all $i\leq n$ we have
    \eqnlab{1.2b3}{
    |\nabla^2 s_i(0)| = 2\Big| \sum_{l=1}^d \lim_{t \rightarrow 0}
    \frac{s_i(te_l)}{t^2} \Big|
    \leq CK_3 \beta i^{-\delta}
    \sum_{l=1}^d \lim_{t \rightarrow 0} \frac{a(te_l)}{t^2}
    =  CK_3 \sigma^2 \beta i^{-\delta}.    }
Note the constant 2 here that is a correction to \cite{HS02}.

Thus, by (\ref{1.2b2}) and Lemma~\ref{lem-cA}
\eqn{|\nabla^2 f_{j}(0)|\leq f_j(0)\sum_{i=1}^j \left[\sigma^2\left(1+CK_2\beta\right)+\frac{CK_3\sigma^2\beta}{i^\delta}\right]\le e^{CK_3\beta} \sigma^2 j
    \Big(1+C(K_2 + K_3 )\beta\Big).}
This completes the proof.
\qed

The next lemma is the key to advancing the induction, as it
provides bounds for $e_{n+1}$ and $g_{n+1}$.

\begin{LEM}
\label{lem-pibds}
 Let $z\in I_{n}$, and assume (H2),
(H3) and (H4). For $k \in [-\pi,\pi]^d$, $2 \leq j \leq n+1$, and
$\epsilon' \in [0,\epsilon]$, the following hold:
\vspace{-2mm}
\begin{tabbing}
(iii) \= \kill (i) \>  $|g_j(k;z)|\leq  K_4' \beta j^{-\theta}$,
\\
(ii) \> $|\nabla^2 g_j(0;z)|\leq   K_4'  \sigma^2 \beta
j^{-\theta+1}$,
\\
(iii) \> $|\partial_z g_j(0;z)|\leq  K_4' \beta j^{-\theta+1},$
\\
(iv) \>  $|g_j(k;z)-g_j(0;z)- a(k) \sigma^{-2}\nabla^2 g_j(0;z)|
\leq  K_4' \beta a(k)^{1+\epsilon'}j^{-\theta+1+\epsilon'},$
\\
(v)   \> $|e_j(k;z)|\leq  K_4' \beta j^{-\theta}$,
\\
(vi)  \> $|e_j(k;z)-e_j(0;z)|\leq  K_4' a(k) \beta
j^{-\theta+1}.$
\end{tabbing}
\end{LEM}

\proof The bounds (\ref{fbdsp}) for $1 \leq m \leq n$
follow from Lemmas~\ref{lem-cA}--\ref{lem-fder}, with
$K=cK_4$ (this defines $c$), assuming that $\beta$ is sufficiently
small. The bounds of the lemma then follow immediately from
Assumptions~\ch{E$_{\theta}$ and G$_{\theta}$}, with $K_4'$ given in (\ref{K4'def}).
\qed

\section{The induction advanced}
\label{sec-adv}
In this section we advance the induction hypotheses (H1--H4) from
$n$ to $n+1$.
Throughout this section, in accordance with the uniformity
condition on (H2--H4), we fix $z \in I_{n+1}$.
We frequently assume $\beta \ll 1$ without explicit comment.

\subsection{Advancement of (H1)}
\label{sec-advH1}

By (\ref{z_n}) and the mean-value theorem,
    \seqn{
    z_{n+1}-z_n & =  -\sum_{m=2}^{n}
    [g_m(0;z_n)-g_m(0;z_{n-1})] -g_{n+1}(0;z_n)\nonumber \\
            & =
    -(z_n-z_{n-1})\sum_{m=2}^{n} \partial_z g_m(0;y_n) -
    g_{n+1}(0;z_{n}),}
for some $y_n$ between $z_n$ and $z_{n-1}$.  By (H1) and
(\ref{Indef}), $y_n \in I_n$. Using Lemma~\ref{lem-pibds} and
(H1), it then follows that
    \seqn{
    |z_{n+1}-z_n| &\leq K_1\beta n^{-\theta}
    \sum\limits_{m=2}^{n}  K_4' \beta m^{-\theta+1} +  K_4' \beta
    (n+1)^{-\theta} \nonumber \\
    & \nonumber \\
    &\leq K_4' \beta(1 + CK_1 \beta) (n+1)^{-\theta}. \label{zdiff}}
Thus (H1) holds for $n+1$, for $\beta$ small and $K_1 > K_4'$.

\medskip
Having advanced (H1) to $n+1$, it then follows from
Lemma~\ref{lem-In} that $I_1 \supset I_2 \supset \cdots \supset
I_{n+1}$.

For $n \geq0$, define
    \eqnlab{zetadef}{
        \zeta_{n+1} =\zeta_{n+1}(z) = \sum_{m=1}^{n+1} g_m(0;z)-1
    = \sum_{m=2}^{n+1} g_m(0;z) + z -1.
    }
The following lemma, whose proof makes use of (H1) for $n+1$, will
be needed in what follows.

\begin{LEM} \label{zetan} For all $z\in I_{n+1}$,
    \seqn{\label{zetanbd}
    |\zeta_{n+1}| & \leq  CK_1 \beta (n+1)^{-\theta+1}.}
\end{LEM}

\proof By (\ref{z_n}) and the mean-value theorem,
    \seqn{
    |\zeta_{n+1}| & =  \Big|(z-z_{n+1}) + \sum_{m=2}^{n+1}
    [g_m(0;z)-g_m(0;z_n)]\Big|\nonumber\\
              & =  \Big|(z-z_{n+1})
    + (z-z_n) \sum_{m=2}^{n+1} \partial_z g_m(0;y_n) \Big|,}
for some $y_n$ between $z$ and $z_n$. Since $z\in I_{n+1} \subset
I_n$ and $z_n \in I_n$, we have $y_n \in I_n$. Therefore, by
Lemma~\ref{lem-pibds},
\eqn{
        |\zeta_{n+1}|  \leq K_1\beta(n+1)^{-\theta+1}
    + K_1\beta n^{-\theta+1}
    \sum_{m=2}^{n+1}  K_4' \beta m^{-\theta+1}
    \leq  K_1\beta(1+CK_4'\beta) (n+1)^{-\theta+1}.
}
The lemma then follows, for $\beta$ sufficiently small.
\qed

\subsection{Advancement of (H2)}
\label{sec-advH1prime}

Let $z\in I_{n+1}$.  As observed in Section~\ref{sec-advH1}, this
implies that $z \in I_j$ for all $j \leq n+1$.
The definitions in (\ref{Delta_n}) imply that
    \eqnlab{vinc}{
    v_{n+1} - v_n = \frac{1}{1+c_{n+1}}(b_{n+1}-b_n) -
    \frac{b_n}{(1+c_n)(1+c_{n+1})}(c_{n+1}-c_n),
    }
with
    \eqnlab{bcdiff}{
    b_{n+1}-b_n = -\frac{1}{\sigma^2}\nabla^2 g_{n+1}(0)
    , \quad
    c_{n+1}-c_n =  n g_{n+1}(0) .
    }
By Lemma~\ref{lem-pibds}, both differences in (\ref{bcdiff}) are bounded by
$K_4' \beta (n+1)^{-\theta+1}$, and, in addition,
\eqnlab{bnear1}{
        |b_j - 1 | \leq C K_4' \beta, \quad
        |c_j | \leq C K_4' \beta
}
for $1 \leq j \leq n+1$.  Therefore
    \eqn{
    |v_{n+1} - v_n|\leq K_2 \beta (n+1)^{-\theta+1},
    }
provided we assume $K_2 \geq 3K_4'$.  This advances (H2).

\subsection{Advancement of (H3)}
\label{sec-advH2}

\subsubsection{The decomposition}

The advancement of the induction hypotheses (H3--H4)
is the most technical part of the proof.
For (H3), we fix $k$ with
$a(k)\leq \gamma (n+1)^{-1}\log{(n+1)}$, and $z\in I_{n+1}$. The
induction step will be achieved as soon as we are able to write
the ratio $f_{n+1}(k)/f_n(k)$ as
\eqn{
        \frac{f_{n+1}(k)}{f_n(k)} = 1-v_{n+1}a(k) + r_{n+1}(k),
}
with $r_{n+1}(0)$ and $r_{n+1}(k)-r_{n+1}(0)$
satisfying the bounds required by (H3).

To begin, we divide the recursion relation (\ref{fkrec}) by
$f_n(k)$, and use (\ref{zetadef}), to obtain
    \seqn{
    \label{rec hat{tau}_n(k)}
    \frac{f_{n+1}(k)}{f_{n}(k)} & =
    1 +\sum_{m=1}^{n+1} \Big[g_m(k)
    \frac{f_{n+1-m}(k)}{f_{n}(k)}-g_m(0)\Big] +\zeta_{n+1}+
    \frac{e_{n+1}(k)}{f_n(k)}.}
By (\ref{Delta_n}),
\eqn{
        v_{n+1} = b_{n+1}-v_{n+1}c_{n+1}
        =-\sigma^{-2} \sum_{m=1}^{n+1}\nabla^2 g_m(0) - v_{n+1}
        \sum_{m=1}^{n+1}(m-1) g_m(0).
        }
Thus we can rewrite (\ref{rec hat{tau}_n(k)}) as
    \eqn{
    \frac{f_{n+1}(k)}{f_{n}(k)} = 1 - v_{n+1} a(k) + r_{n+1}(k),
    \label{the eq}
    }
where
    \eqn{
    r_{n+1}(k) = X(k)+Y(k)+Z(k)+ \zeta_{n+1}
    }
with
    \seqn{
    X(k)    & = \sum_{m=2}^{n+1}
    \Big[\big(g_m(k)-g_m(0)\big)\frac{f_{n+1-m}(k)}{f_n(k)}-a(k)
    \sigma^{-2}\nabla^2g_m(0)\Big],\\
    \label{H2IIdef}
    Y(k) & = \sum_{m=2}^{n+1} g_m(0)\left[
    \frac{f_{n+1-m}(k)}{f_n(k)}-1- (m-1) v_{n+1}a(k)\right],
    \hspace{6mm} \\
    Z(k) & =  \frac{e_{n+1}(k)}{f_n(k)}.}
The $m=1$ terms in $X$ and $Y$ vanish and have not been
included.

We will prove that
\eqnlab{rbds}{
    |r_{n+1}(0)|\leq \frac{C(K_1+ K_4') \beta}{(n+1)^{\theta-1}},
    \qquad \qquad
    |r_{n+1}(k)-r_{n+1}(0)|
    \leq \frac{C K_4' \beta a(k)}{(n+1)^{\delta}}.
}
This gives (H3) for $n+1$, provided we assume that $K_3 \gg K_1$
and $K_3 \gg K_4'$.
To prove the bounds on $r_{n+1}$ of (\ref{rbds}), it will be convenient
to make use of some elementary convolution bounds, as well as some bounds
on ratios involving $f_j$.  These preliminary bounds are given
in Section~\ref{sec-ratiobds}, before we present the proof of
(\ref{rbds}) in Section~\ref{sec-XYZ}.

\subsubsection{Convolution and ratio bounds}
\label{sec-ratiobds}

The proof of (\ref{rbds}) will make use of the following
elementary convolution bounds.  To keep the discussion simple, we
do not obtain optimal bounds.

\begin{LEM}
    \label{lem-conv}
For $n \geq 2$,
    \eqnlab{conv-bound}{
    \sum_{m=2}^{n} \frac{1}{m^a}
    \sum_{j=n-m+1}^n \frac{1}{j^b} \leq \left\{
    \begin{array}{lll}&C n^{-(a\wedge b)+1} &\mbox{for }a,b>1\\
    &Cn^{-(a-2)\wedge b}    &\mbox{for }a>2, b>0\\
    &Cn^{-(a-1)\wedge b}    &\mbox{for }a>2, b>1\\
    &Cn^{-a\wedge b}    &\mbox{for }a,b>2.\end{array}\right.
    }
\end{LEM}

\proof
Since $m+j\geq n$, either $m$ or $j$ is at least $\frac{n}{2}$.
Therefore
    \eqnlab{conv-bound2}{
    \sum_{m=2}^{n} \frac{1}{m^a}
    \sum_{j=n-m+1}^n \frac{1}{j^b} \leq \left(\frac{2}{n}\right)^a
    \sum_{m=2}^{n}\sum_{j=n-m+1}^n \frac{1}{j^b}
    +\left(\frac{2}{n}\right)^b \sum_{m=2}^{n}\sum_{j=n-m+1}^n
    \frac{1}{m^a}.
    }
If $a,b>1$, then the first term is bounded by $Cn^{1-a}$ and the
second by $Cn^{1-b}$.\\ If $a>2, b>0$, then the first term is
bounded by $Cn^{2-a}$ and the second by $Cn^{-b}$.\\ If $a>2,
b>1$, then the first term is bounded by $Cn^{1-a}$ and the second
by $Cn^{-b}$.\\ If $a,b>2$, then the first term is bounded by
$Cn^{-a}$ and the second by $Cn^{-b}$.
\qed

We also will make use of several estimates involving ratios.
We begin with some preparation.
Given a vector $x=(x_l)$ with $\sup_l|x_l|<1$, define $\chi(x) = \sum_l
\frac{|x_l|}{1-|x_l|}$.
The bound $(1-t)^{-1} \leq \exp [t(1-t)^{-1}]$, together with Taylor's
Theorem applied to $f(t) = \prod_l \frac{1}{1-tx_l}$, gives
\eqnlab{Taylor1}{
        \left|\prod_l \frac{1}{1-x_l} -1\right| \leq \chi(x) e^{\chi(x)},
        \quad
        \left|\prod_l \frac{1}{1-x_l} -1 -\sum_l x_l\right|
        \leq \chi(x)^2 e^{\chi(x)}
}
as follows.
Firstly,
\eqn{\frac{df}{dt}=f(t)\sum_{j=1}^d \frac{x_j}{1-tx_j} =\left[\prod_{l=1}^d\R{1-tx_l}\right]\sum_{j=1}^d \frac{x_j}{1-tx_j}\le \left[\prod_{l=1}^d e^{\frac{|tx_j|}{1-|tx_j|}}\right]\sum_{j=1}^d \frac{|x_j|}{1-|tx_j|},}
which gives $f'(0)=\sum_{j=1}^d x_j$, and for $|t|\le 1$, $|f'(t)|\le \chi(x)e^{\chi(x)}$.  This gives the first bound by Taylor's Theorem.  The second bound can be obtained in the same way using the fact that
\eqn{\frac{d^2f}{dt^2}=f(t)\left[\sum_{j=1}^d\frac{x_j^2}{(1-tx_j)^2}+ \left(\sum_{j=1}^d\frac{x_j}{1-tx_j}\right)^2\right].}
We assume throughout the rest of this section that $a(k) \leq \gamma
(n+1)^{-1} \log (n+1)$ and $2 \leq m \leq n+1$, and define
\eqnlab{chidef}{
        \psi_{m,n} = \sum_{j=n+2-m}^n \frac{|r_j(0)|}{1-|r_j(0)|}, \quad
        \chi_{m,n}(k) = \sum_{j=n+2-m}^n \
        \frac{v_j a(k)+ |s_j(k)|}{1-v_j a(k) - |s_j(k)|}.
}
By (\ref{vnear1}) and (\ref{sbd}),
    \eqnlab{chibd1}{\chi_{m,n}(k) \leq (m-1) a(k)Q(k)
    \quad \mbox{ with } \quad Q(k) = [1+C(K_2+K_3)\beta][1+Ca(k)],}
where we have used the fact that for $|x|\le \hlf$, $\R{1-x}\le 1+2|x|$.  In our case $x=v_j a(k)+ |s_j(k)|$ satisfies $|x|\le (1+CK_2\beta)a(k)+CK_3\beta a(k)$.
Since $a(k) \leq \gamma (n+1)^{-1} \log (n+1)$, we have
$Q(k) \leq [1+C(K_2+K_3)\beta][1+C\gamma (n+1)^{-1} \log (n+1)]$.
Therefore
\seqnlab{chibd3}{
        e^{\chi_{m,n}(k)} &\leq e^{\gamma \log(n+1) Q(k)}
        \leq e^{\gamma \log(n+1)[1+C(K_2+K_3)\beta]}e^{\frac{C\gamma^2(\log (n+1))^2}{n+1}}\\
&\le e^{\gamma \log(n+1)[1+C(K_2+K_3)\beta]}e^{4C\gamma^2}\le
        C(n+1)^{\gamma q},
}
where we have used the fact that $\log x \le 2\sqrt{x}$, and where $q=1+C(K_2+K_3)\beta$ may be taken to be as close to $1$ as desired,
by taking $\beta$ to be small.

We now turn to the ratio bounds.  It follows from (H3)
and the first inequality of (\ref{Taylor1}) that
\seqnlab{ratio1.a}{
        \left| \frac{f_{n+1-m}(0)}{f_n(0)} - 1 \right|&=\left|\prod_{i=n+2-m}^n\R{1-(-r_i(0))}-1\right|\\
        &\leq \psi_{m,n}e^{\psi_{m,n}}
        \leq
        \sum_{j=n+2-m}^n \frac{CK_3 \beta}{j^{\theta-1}}
        \leq
        \frac{CK_3 \beta}{(n+2-m)^{\theta-2}}}
Therefore
\eqnlab{ratio0}{
        \left| \frac{f_{n+1-m}(0)}{f_n(0)} \right|
        \leq
        1+CK_3 \beta.
}
By (\ref{fs}),
\seqn{
        \left| \frac{f_{n+1-m}(k)}{f_n(k)} - 1 \right|&=\left|\frac{f_{n+1-m}(0)}{f_n(0)}\prod_{j=n+2-m}^n \R{[1-v_ja(k) +s_j(k)]}-\frac{f_{n+1-m}(0)}{f_n(0)}+\frac{f_{n+1-m}(0)}{f_n(0)}-1\right|\\
        &\leq
        \left| \frac{f_{n+1-m}(0)}{f_n(0)} \right|
        \left| \prod_{j=n+2-m}^n \R{[1-v_ja(k) +s_j(k)]}
        - 1 \right|
        + \left| \frac{f_{n+1-m}(0)}{f_n(0)} - 1 \right|.
}
The first inequality of (\ref{Taylor1}), together with
(\ref{chibd1}--\ref{ratio0}), then gives
\eqnlab{ratio1}{
        \left| \frac{f_{n+1-m}(k)}{f_n(k)} - 1 \right|
        \leq
        C (m-1) a(k)(n+1)^{\gamma q} + \frac{CK_3 \beta}{(n+2-m)^{\theta-2}} .
}
Similarly,
\eqnlab{ratio2}{
        \left| \frac{f_{n}(0)}{f_n(k)} - 1 \right|=\left|\prod_{i=1}^n\R{1-v_ja(k) +s_j(k)}-1\right|
        \leq
        \chi_{n+1,n}(k) e^{\chi_{n+1,n}(k)}
        \leq
        C a(k) (n+1)^{1+\gamma q}.
}

Next, we estimate the quantity $R_{m,n}(k)$, which is defined by
\eqnlab{Rdef}{
        R_{m,n}(k) = \prod_{j=n+2-m}^n [1-v_ja(k) +s_j(k)]^{-1}
        - 1
        - \sum_{j=n+2-m}^n [v_ja(k) -s_j(k)] .
}
By the second inequality of (\ref{Taylor1}), together with
(\ref{chibd1}) and (\ref{chibd3}), this obeys
\eqnlab{Rbd}{
        |R_{m,n}(k)| \leq  \chi_{m,n}(k)^2 e^{\chi_{m,n}(k)}
        \leq
        C m^2 a(k)^2 (n+1)^{\gamma q}.
}

Finally, we apply (H3) with $\R{1-x}-1=\frac{x}{1-x}\le \frac{|x|}{1-|x|}$ to obtain for $m\le n$,
\eqnlab{ratio3}{
        \left| \frac{f_{m-1}(k)}{f_m(k)} - 1 \right|
        =
        \left| [1-v_ma(k) +(r_m(k)-r_m(0)) + r_m(0)]^{-1}
        - 1 \right|
        \leq
        Ca(k) + \frac{CK_3 \beta}{m^{\theta-1}} .
}
Note that for example, $1-(|v_ma(k)| +|r_m(k)-r_m(0)| + |r_m(0)|)>c$ for small enough $\beta$ (depending on $\gamma$, among other things).

\subsubsection{The induction step}
\label{sec-XYZ}

By definition,
    \eqnlab{rp0}{
    r_{n+1}(0) = Y(0)+Z(0)+\zeta_{n+1}
    }
and
    \eqnlab{rpk0}{
    r_{n+1}(k)-r_{n+1}(0) = X(k) + \Big( Y(k)-Y(0) \Big)
     + \Big( Z(k)-Z(0) \Big).
    }
Since $|\zeta_{n+1}| \leq CK_1 \beta (n+1)^{-\theta+1}$ by
Lemma~\ref{zetan}, to prove (\ref{rbds}) it suffices to show that
    \eqnlab{rp0suf}{
    |Y(0)| \leq C K_4' \beta (n+1)^{-\theta+1} ,\quad
    |Z(0)| \leq C K_4' \beta (n+1)^{-\theta+1}
    }
and
    \seqn{
    \label{rpk0suf}
    & |X(k)|  \leq  C K_4' \beta a(k) (n+1)^{-\delta}, \quad
    |Y(k)-Y(0)| \leq C K_4' \beta a(k) (n+1)^{-\delta},
    \nonumber \\ & \hspace{25mm}
    |Z(k)-Z(0)|  \leq  C K_4' \beta a(k) (n+1)^{-\delta}.}
The remainder of the proof is devoted to establishing
(\ref{rp0suf}) and (\ref{rpk0suf}).

\medskip\noindent
{\em Bound on $X$}. We write $X$ as $X=X_1+X_2$, with
    \seqn{
    \label{X1def}
    X_1 &= \sum_{m=2}^{n+1}
    \Big[g_m(k)-g_m(0)-a(k) \sigma^{-2}\nabla^2g_m(0)\Big] ,
    \\
    X_2 &=\sum_{m=2}^{n+1}
    \Big[g_m(k)-g_m(0)\Big]\Big[\frac{f_{n+1-m}(k)}{f_n(k)}-1\Big].}
The term $X_1$ is bounded using Lemma~\ref{lem-pibds}(iv)
with $\epsilon' \in (\delta, \epsilon )$, and
using the
fact that $a(k)\leq \gamma (n+1)^{-1} \log{(n+1)}$, so that $a(k)^{\epsilon'}\le \left(\frac{\gamma \log (n+1)}{n+1}\right)^{\epsilon'}\le \sm{C}{(n+1)}{\delta}$ by
\eqn{
    |X_1|  \leq  K_4' \beta
    a(k)^{1+\epsilon'}\sum_{m=2}^{n+1}
    \frac{1}{m^{\theta-1-\epsilon'}} \leq C K_4' \beta
    a(k)^{1+\epsilon'} \leq \frac{CK_4' \beta a(k)}{(n+1)^{\delta}}.
\label{I}
}

For $X_2$, we first apply Lemma~\ref{lem-pibds}(ii,iv), with
$\epsilon' =0$, to obtain
\eqn{
        |g_m(k)-g_m(0)| \leq 2 K_4' \beta a(k) m^{-\theta+1}.
}
Applying (\ref{ratio1}) then gives
\eqnlab{X2bd}{
        |X_2| \leq CK_4' \beta a(k)
        \sum_{m=2}^{n+1} \frac{1}{m^{\theta-1}}
        \left(
        (m-1) a(k)(n+1)^{\gamma q}
        + \frac{K_3\beta}{(n+2-m)^{\theta-2}}
        \right).
}
By the elementary estimate
\eqn{\sum_{m=2}^{n+1}\sm{1}{m}{\theta-1}\sm{1}{(n+2-m)}{\theta-2}\le \sm{C}{(n+1)}{\theta-2},}
which is proved easily by breaking the sum up according to $m \le \lfloor \frac{n+1}{2}\rfloor$, the contribution from the
second term on the right side
is bounded above by $CK_3K_4'\beta^2 a(k) (n+1)^{-\theta+2}$.
The first term is bounded above by
\eqn{
        CK_4' \beta a(k) (n+1)^{\gamma q - 1} \log (n+1)
        \times \begin{cases} (n+1)^{0 \vee (3-\theta)} & (\theta \neq 3)
        \\ \log(n+1) & (\theta=3) . \end{cases}
}
Since we may choose $q$ to be as close to $1$ as desired,
and since $\delta + \gamma < 1 \wedge (\theta-2)$ by (\ref{agddef}),
this is bounded above by $CK_4' \beta a(k)(n+1)^{-\delta}$.
With (\ref{I}), this proves the bound on $X$ in (\ref{rpk0suf}).

\medskip \noindent {\em Bound on $Y$}.
By (\ref{fs}),
    \eqnlab{AkA0ratio}{
    \frac{f_{n+1-m}(k)}{f_n(k)} =
    \frac{f_{n+1-m}(0)}{f_n(0)} \prod_{j=n+2-m}^n
    [1 -v_j a(k)+s_j(k)]^{-1}  .
    }
Recalling the definition of $R_{m,n}(k)$
in (\ref{Rdef}), we can therefore decompose $Y$ as $Y=Y_1+Y_2+Y_3+Y_4$ with
    \seqn{
    Y_1 & =   \sum_{m=2}^{n+1} g_m(0)
    \frac{f_{n+1-m}(0)}{f_n(0)}
    R_{m,n}(k),  \\
    Y_2 & =   \sum_{m=2}^{n+1} g_m(0)
    \frac{f_{n+1-m}(0)}{f_n(0)}
    \sum_{j=n+2-m}^n \left[ (v_j-v_{n+1})a(k)
    - s_j(k) \right], \\
    Y_3 & =   \sum_{m=2}^{n+1} g_m(0)
    \left[ \frac{f_{n+1-m}(0)}{f_n(0)} - 1 \right]
    (m-1)v_{n+1} a(k), \\
    Y_4 & =   \sum_{m=2}^{n+1} g_m(0)
    \left[ \frac{f_{n+1-m}(0)}{f_n(0)} - 1 \right].}
Then
    \eqn{
    Y(0) = Y_4  \quad \mbox{ and } \quad
    Y(k)-Y(0) = Y_1+Y_2+Y_3.
    }

\noindent For $Y_1$, we use Lemma~\ref{lem-pibds}, (\ref{ratio0}) and (\ref{Rbd})
to obtain
    \eqnlab{III2z}{
    |Y_1| \leq  CK_4'\beta  a(k)^2 (n+1)^{\gamma q}
    \sum\limits_{m=2}^{n+1} \frac{1}{m^{\theta-2}}.
    }
As in the analysis of the first term of (\ref{X2bd}), we therefore have
    \eqnlab{III2a}{
        |Y_1| \leq  \frac{CK_4' \beta a(k)}{(n+1)^{\delta}}.
    }
For $Y_2$, we use $\theta-2>\delta>0$ with Lemma~\ref{lem-pibds}, (\ref{ratio0}), (H2) (now established up to $n+1$), (\ref{sbd}) and Lemma~\ref{lem-conv}
to obtain
    \eqnlab{II2bd}{
    |Y_2| \leq \sum_{m=2}^{n+1} \frac{ K_4' \beta}{m^{\theta}}
    C
    \sum_{j=n+2-m}^n \left[ \frac{K_2 \beta a(k)}{j^{\theta-2}}
    + \frac{K_3 \beta a(k)}{j^{\delta}} \right]
    \leq \frac{CK_4' (K_2+K_3) \beta^2 a(k)}{(n+1)^{\delta}}.
    }
The term $Y_3$ obeys
    \eqnlab{II3bd}{
    |Y_3| \leq
    \sum_{m=2}^{n+1} \frac{K_4'\beta}{m^{\theta-1}}
    \frac{CK_3 \beta }{(n+2-m)^{\theta-2}}  a(k)
    \leq \frac{CK_4' K_3 \beta^2 a(k)}{(n+1)^{\theta-2}},
    }
where we used Lemma~\ref{lem-pibds}, (\ref{ratio1.a}),
(\ref{vnear1}), and an elementary convolution bound.   This proves
the bound on $|Y(k)-Y(0)|$ of (\ref{rpk0suf}), if $\beta$ is
sufficiently small.

We bound $Y_4$ in a similar fashion, using Lemma~\ref{lem-conv}
and the intermediate bound of (\ref{ratio1.a}) to obtain
    \eqn{
    \label{Vest}
    |Y_4|  \leq
    \sum_{m=2}^{n+1} \frac{ K_4'\beta}{m^{\theta}}
    \sum_{j = n+2-m}^n
    \frac{CK_3 \beta}{j^{\theta-1}}
    \leq \frac{CK_4'K_3 \beta^2 }{(n+1)^{\theta-1}}.
    }
Taking $\beta$ small then gives the bound on $Y(0)$ of
(\ref{rp0suf}).

\medskip \noindent {\em Bound on $Z$}.
We decompose $Z$ as
\eqnlab{IIIsplit}{
    Z = \frac{e_{n+1}(0)}{f_n(0)}
    +\frac{1}{f_n(0)}
    \left[ e_{n+1}(k) - e_{n+1}(0) \right]
    + \frac{e_{n+1}(k)}{f_n(0)}\left[\frac{f_n(0)}{f_n(k)} - 1\right]
    = Z_1+Z_2+Z_3.
}
Then
        \eqn{
        Z(0)=Z_1 \quad \mbox{ and } \quad Z(k)-Z(0) = Z_2+Z_3.
        }
Using Lemma~\ref{lem-pibds}(v,vi), and (\ref{ratio0})
with $m=n+1$, we
obtain
    \eqn{
    |Z_1| \leq CK_4' \beta (n+1)^{-\theta} \quad \mbox{ and }
    \quad
    |Z_2| \leq CK_4' \beta a(k)(n+1)^{-\theta+1}.
    }
Also, by Lemma~\ref{lem-pibds}, (\ref{ratio0}) and (\ref{ratio2}), we have
\eqn{
     |Z_3|
     \leq   CK_4' \beta (n+1)^{-\theta}
     a(k) (n+1)^{1 + \gamma q}
     \leq CK_4' \beta a(k)(n+1)^{-(1+\delta)},
}
for small enough $q$, where we again use $\gamma+\delta<\theta-2$.
%

This completes the proof of (\ref{rbds}), and hence completes
the advancement of (H3) to $n+1$.

\subsection{Advancement of (H4)}
\label{sec-advH3}\label{sec-advH35}

In this section, we fix $a(k) > \gamma (n+1)^{-1} \log (n+1)$. To
advance (H4) to $j=n+1$, we first recall the definitions of
$b_{n+1}$, $\zeta_{n+1}$ and $X_1$ from (\ref{Delta_n}),
(\ref{zetadef}) and (\ref{X1def}).  After some algebra,
(\ref{fkrec}) can be rewritten as
    \eqn{
    f_{n+1}(k) =  f_n(k)\Big( 1~- a(k)b_{n+1}
    + X_1 + \zeta_{n+1}\Big) + W +e_{n+1}(k),
    \label{the eq (H3-H4)}
    }
with
    \eqn{
%
    \label{H3IIdef}
    W   = \sum_{m=2}^{n+1} g_m(k) \left[f_{n+1-m}(k)-f_n(k)\right].
    }

We already have estimates for most of the relevant terms. By
Lemma~\ref{zetan}, we have $|\zeta_{n+1}| \leq CK_1
\beta(n+1)^{-\theta+1}$. By (\ref{I}), $|X_1| \leq CK_4' \beta
a(k)^{1+\epsilon'}$, for any $\epsilon' \in(\delta, \epsilon )$. By
Lemma~\ref{lem-pibds}(v), $|e_{n+1}(k)|\leq K_4' \beta
(n+1)^{-\theta}$. It remains to estimate $W$.  We will show below that $W$
obeys the bound
    \eqn{
    \label{Wbd}
    |W|
    \leq \frac{CK_4' \beta}{a(k)^{a-1}(n+1)^{\theta}}
    (1+K_3\beta +K_5).
    }
Before proving (\ref{Wbd}), we will first show that it is sufficient
for the advancement of (H4).

In preparation for this, we first note that it suffices to consider only
large $n$.   In fact, since
$|f_n(k;z)|$ is bounded uniformly in $k$ and in $z$ in a
compact set by Assumption~S, and since $a(k)\leq 2$, it is clear
that both inequalities of (H4) hold for all $n\leq N$, if we
choose $K_4$ and $K_5$ large enough (depending on $N$).
We therefore assume in the following
that $n \geq N$ with $N$ large.

Also, care is required to invoke (H3) or (H4), as applicable, in
estimating the factor $f_n(k)$ of (\ref{the eq (H3-H4)}).  Given $k$,
(H3) should be used for the value $n$ for which
$\gamma (n+1)^{-1} \log (n+1) < a(k) \leq \gamma n^{-1} \log n$
((H4) should be used for larger $n$).
We will now show that
the bound of (H3) actually implies the first
bound of (H4) in this case.  To see this, we use Lemma~\ref{lem-cA}
to see that there are $q,q'$ arbitrarily close to $1$ such that
\eqn{
\label{H3toH4}
        |f_n(k)| \leq Ce^{-qa(k)n} \leq \frac{C}{(n+1)^{q\gamma n/(n+1)}}
        \leq \frac{C}{n^{q' \gamma}}
        \leq \frac{C}{n^{\theta}}\frac{n^{\lambda}}{n^{q'\gamma +\lambda -\theta}}
        \leq \frac{C}{n^{\frac{d}{2p}}a(k)^{\lambda}},
}
where we used the fact that $\gamma + \lambda -\theta >0$ by (\ref{agddef}).
Thus, taking $K_4 \gg 1$,
we may use the
first bound of (H4) also for the value of $n$ to which
(H3) nominally applies.  We will do so in what follows, without further comment.

\medskip\noindent {\em Advancement of the second bound of (H4)
assuming (\ref{Wbd})}.
To advance the second estimate in (H4),
we use (\ref{the eq (H3-H4)}), (H4), and the bounds found
above, to obtain
    \seqn{
    \Big|f_{n+1}(k) -  f_n(k)\Big|
    &\leq  \big|f_n(k)\big|~
    \big|-a(k)b_{n+1}+ X_1 + \zeta_{n+1}\big| + |W|+|e_{n+1}(k)| \nonumber \\
    & \leq \frac{K_4}{n^{\theta}a(k)^{\lambda}}
    \left(a(k)b_{n+1}+ C K_4' \beta a(k)^{1+\epsilon'}
    + \frac{CK_1\beta }{(n+1)^{\theta-1}}
    \right)\nn\\
    \label{H3sec}
    &\quad+ \frac{CK_4' \beta (1+K_3\beta +K_5)}{(n+1)^{\theta}a(k)^{\lambda-1}}
    + \frac{ K_4' \beta}{(n+1)^{\theta}}.}
Since $b_{n+1}=1+\mc{O}(\beta)$ by (\ref{bnear1}), and
since $(n+1)^{-\theta+1} < [a(k) /\gamma\log (n+1)]^{\theta-1} \leq C a(k)$, the
second estimate in (H4) follows for $n+1$ provided $K_5 \gg K_4$
and $\beta$ is sufficiently small.

\medskip\noindent {\em Advancement of the first bound of (H4) assuming
(\ref{Wbd}}).
To advance the first estimate of (H4), we argue as in
(\ref{H3sec}) to obtain
    \seqn{
    \big|f_{n+1}(k)\big|  &\leq
    \big|f_n(k)\big|~\Big|1~-a(k)b_{n+1}+ X_1 +\zeta_{n+1}\Big|
    + |W| +|e_{n+1}(k)| \nonumber\\
    &\leq  \frac{K_4}{n^{\theta}a(k)^{\lambda}}
    \left( |1~-a(k)b_{n+1}|
    + C K_4' \beta a(k)^{1+\epsilon'}
    + \frac{CK_1\beta }{(n+1)^{\theta-1}}
    \right)\nn\\
    &\quad+ \frac{CK_4' \beta (1+K_3\beta +K_5)}{(n+1)^{\theta}a(k)^{\lambda-1}}
    + \frac{ K_4' \beta}{(n+1)^{\theta}}.
    \label{H3bd1}}
We need to argue that the right-hand side is no larger than $K_4
(n+1)^{-\theta} a(k)^{-\lambda}$. To achieve this, we will use
separate arguments for $a(k) \leq \frac 12$ and $a(k) > \frac 12$.
These arguments will be valid only when $n$ is large enough.

Suppose that $a(k) \leq \frac 12$. Since $b_{n+1} = 1+
\mc{O}(\beta)$ by (\ref{bnear1}), for $\beta$ sufficiently small we
have
    \eqn{
    1~-b_{n+1}a(k) \geq 0.
    }
Hence, the absolute value signs on the right side of (\ref{H3bd1})
may be removed.
Therefore, to obtain the first estimate of (H4) for $n+1$,
it now suffices to show that
    \eqn{
    \label{H3bd2}
    1~-ca(k)
    +\frac{CK_1\beta}{(n+1)^{\theta-1}}
    \leq \frac{n^{\theta}}{(n+1)^{\theta}},
    }
for $c$ within order $\beta$ of 1.  The term $ca(k)$ has been
introduced to absorb $b_{n+1}a(k)$, the order $\beta$ term in
(\ref{H3bd1}) involving $a(k)^{1+\epsilon'}$, and the last two
terms of (\ref{H3bd1}).  However,
$a(k) > \gamma (n+1)^{-1} \log(n+1)$.  From this, it can be
seen that (\ref{H3bd2})
holds for $n$ sufficiently large and $\beta$ sufficiently small.

Suppose, on the other hand, that $a(k) > \frac 12$. By
(\ref{Dbound3}), there is a positive $\eta$, which we
may assume lies in $(0,\frac{1}{2})$, such that $-1+\eta <1-a(k)<
\frac 12$. Therefore $|1-a(k)| \leq 1-\eta$ and
    \eqn{
    |1~-b_{n+1}a(k)|
    \leq |1-a(k)| + |b_{n+1}-1| \, |a(k)|
    \leq 1-\eta
    + 2|b_{n+1}-1|.
    }
Hence
    \eqn{
    |1-a(k)b_{n+1}| + C K_4' \beta a(k)^{1+\epsilon'}
    +\frac{CK_1\beta }{(n+1)^{\theta-1}}
    \leq 1-\eta + C(K_1+K_4')\beta,
    }
and the right side of (\ref{H3bd1}) is at most
    \seqn{
    &\frac{K_4}{n^{\theta}a(k)^{\lambda}}
    \left[1-\eta + C(K_1 +K_4')\beta \right] +
    \frac{CK_4'(1+K_3\beta + K_5) \beta}{(n+1)^{\theta}a(k)^{\lambda}}\nonumber\\
    &\quad\leq  \frac{K_4}{n^{\theta}a(k)^{\lambda}}
    \left[1-\eta + C(K_5K_4'+K_1)\beta \right].}
This is less than $K_4 (n+1)^{-\theta}a(k)^{-\lambda}$ if $n$ is
large and $\beta$ is sufficiently small.

\noindent This advances the first bound in (H4), assuming (\ref{Wbd}).

\medskip\noindent {\em Bound on $W$}.
We now obtain the bound (\ref{Wbd}) on $W$.  As a first step,
we rewrite $W$ as
    \eqn{
    \label{Wdef}
    W  =\sum_{j=0}^{n-1} g_{n+1-j}(k)\sum_{l=j+1}^n
    [f_{l-1}(k)-f_l(k)].
    }
Let
    \eqn{
    m(k) = \begin{cases}
    1 & (a(k) > \gamma 3^{-1}\log 3) \\
    \max \{ l\in \{3,\ldots, n\} : a(k) \leq \gamma l^{-1} \log l \}
    & ( a(k) \leq \gamma 3^{-1}\log 3).
    \end{cases}
    }

For $l \leq m(k)$, $f_l$ is in the domain of (H3), while for $l>m(k)$,
$f_l$ is
in the domain of (H4).  By hypothesis, $a(k) > \gamma (n+1)^{-1} \log(n+1)$.
We divide the sum over $l$ into two parts,
corresponding respectively to $l \leq m(k)$ and $l> m(k)$,
yielding $W=W_1 + W_2$.
By Lemma~\ref{lem-pibds}(i),
    \seqn{
    |W_1|  & \leq \sum_{j=0}^{m(k)} \frac{ K_4' \beta}{(n+1-j)^{\theta}}
    \sum_{l=j+1}^{m(k)} |f_{l-1}(k)-f_l(k)|
    \\
    |W_2|  & \leq \sum_{j=0}^{n-1}\frac{ K_4' \beta}{(n+1-j)^{\theta}}
    \sum_{l=(m(k) \vee j)+1}^n |f_{l-1}(k)-f_l(k)|.}
The term $W_2$ is easy, since by (H4) and Lemma \ref{lem-conv} we have
    \eqn{
    \label{sum1}
    |W_2| \leq \sum\limits_{j=0}^{n-1} \frac{
    K_4' \beta}{(n+1-j)^{\theta}} \sum\limits_{l=j+1}^n
    \frac{K_5}{a(k)^{\lambda-1} \; l^{\theta}} \leq \frac{C K_5 K_4'
    \beta}{ a(k)^{\lambda-1} (n+1)^{\theta}}.
    }



\noindent For $W_{1}$, we have the estimate
\eqn{
\label{W1pbd}
    |W_1| \leq \sum_{j=0}^{m(k)} \frac{ K_4' \beta}{(n+1-j)^{\theta}}
    \sum_{l=j+1}^{m(k)} |f_{l-1}(k)-f_l(k)|.
}
For $1 \leq l \leq m(k)$, it follows from Lemma~\ref{lem-cA}
and (\ref{ratio3}) that
    \eqn{
    \label{fldiff}
    |f_{l-1}(k)-f_l(k)|
    \leq
    C e^{- q a(k) l}\left( a(k) + \frac{K_3 \beta}{l^{\theta-1}}\right),
    }
with $q =1-\mc{O}(\beta)$.
We fix a small $r > 0$, and bound the summation over $j$
in (\ref{W1pbd}) by summing separately over $j$ in the ranges $0
\leq j \leq (1-r)n$ and $(1-r)n \leq j \leq m(k)$ (the
latter range may be empty).  We denote the contributions from
these two sums by $W_{1,1}$ and $W_{1,2}$ respectively.

To estimate $W_{1,1}$, we will make use of the bound
\eqn{
        \sum_{l=j+1}^\infty e^{- q a(k) l}  l^{-b}
        \leq
        C e^{- q a(k) j} \quad \quad (b>1).
}
With (\ref{W1pbd}) and (\ref{fldiff}), this gives
    \seqn{
    |W_{1,1}|
    & \leq  \frac{CK_4' \beta}{(n+1)^{\theta}}
    \sum_{j=0}^{(1-r)n}
    e^{- q a(k) j}
    \left( 1 + K_3 \beta\right)
    \nn \\
    \label{W11p}
    & \leq  \frac{CK_4' \beta}{(n+1)^{\theta}} \frac{1+K_3\beta}{a(k)}
    \leq \frac{CK_4' \beta}{(n+1)^{\theta}} \frac{1+K_3\beta}{a(k)^{\lambda-1}}.}

For $W_{1,2}$, we have
    \eqn{
    |W_{1,2}| \leq   \sum_{j=(1-r)n}^{m(k)}
    \frac{CK_4'\beta }{(n+1-j)^{\theta}}
    \sum_{l=j+1}^{m(k)}
     e^{-q a(k) l}\left(a(k) + \frac{K_3 \beta}{l^{\theta-1}}\right).
    }
Since $l$ and $m(k)$ are comparable ($(1-r)(n+1)<(1-r)n+1\le l \le m(k)<n+1$) and large, it follows as in
(\ref{H3toH4}) that 
\eqn{
    e^{-q a(k) l}
    \left(a(k) + \frac{K_3 \beta}{l^{\theta-1}}\right)
    \leq \frac{C}{a(k)^{\lambda} l^{\theta}}
    \left(a(k) + \frac{K_3 \beta}{l^{\theta-1}}\right)
    \leq \frac{C(1+K_3\beta)}{a(k)^{\lambda-1}l^{\theta}},
}
where we have used the definition of $m(k)$ in the form $\frac{\gamma \log (m(k)+1)}{m(k)+1}<a(k)\le \frac{\gamma \log (m(k))}{m(k)}$ as well as the facts that $\lambda>\theta-\gamma$ and that $q(1-r)$ can be chosen as close to 1 as we like to obtain the intermediate inequality, and the same bound on $a(k)$ together with the fact that $\theta>2$ to obtain the last inequality.
Hence, by Lemma~\ref{lem-conv},
    \eqn{
    \label{W12p}
    |W_{1,2}| \leq \frac{C(1+K_3\beta)K_4' \beta}{a(k)^{\lambda-1}}
    \sum_{j=(1-r)n}^{m(k)} \frac{1}{(n+1-j)^{\theta}}
    \sum_{l=j+1}^{m(k)} \frac{1}{l^{\theta}}
    \leq \frac{C(1+K_3\beta)K_4' \beta}{a(k)^{\lambda-1}(n+1)^{\theta}}.
    }

Summarising, by (\ref{W11p}), (\ref{W12p}),
and (\ref{sum1}), we have
    \eqn{
    |W| \leq |W_{1,1}| + |W_{1,2}|  +|W_2|
    \leq \frac{CK_4' \beta}{a(k)^{\lambda-1}(n+1)^{\theta}}
    (1+K_3\beta +K_5),
    }
which proves (\ref{Wbd}).

\section{Proof of the main results}
\label{sec-pf}

As a consequence of the completed induction, it follows from
Lemma~\ref{lem-In} that $I_1 \supset I_2 \supset I_3 \supset
\cdots$, so $\cap_{n=1}^\infty I_n$ consists of a single point
$z=z_c$.  Since $z_0=1$, it follows from (H1) that $z_c=1+\mc{O}(\beta)$.
We fix $z=z_c$ throughout this section.  The constant $A$ is
defined by $A = \prod_{i=1}^\infty [1+
r_i(0)] = 1+ {\cal O}(\beta)$.
By (H2), the sequence $v_n(z_c)$ is a Cauchy sequence.
The constant $v$ is defined to be the limit of this Cauchy sequence.
By (H2), $v = 1+ {\cal O}(\beta)$ and
        \eqn{
        \label{Delta_nlimit}
        |v_n(z_c) - v| \leq {\cal O}( \beta n^{-\theta+2}).
        }

\subsection{Proof of Theorem \ref{thm-1p}}

\smallskip \noindent {\em Proof of Theorem \ref{thm-1p}(a).}
By (H3),
        \eqn{
        \label{hat{tau}_nlimit}
        |f_n(0;z_c) - A|
        = \prod_{i=1}^n [1+r_i(0)]
        \big| 1 - \prod_{i=n+1}^\infty [1+r_i(0)]\big|
        \leq {\cal O}(\beta n^{-\theta+2}).
        }
Suppose $k$ is such that $a(k/\sqrt{\sigma^2 vn}) \leq \gamma n^{-1} \log
n$, so that (H3) applies.  Here, we use the $\gamma$ of (\ref{agddef}).
By (\ref{momentD}),
$a(k) = \sigma^2k^2/2d  + {\cal O}(k^{2+2\epsilon})$ with $\epsilon > \delta$,
where we now allow constants in error terms to depend on $L$.
Using this, together with (\ref{fs}--\ref{sbd}),
\ref{Delta_nlimit}, and  $\delta <
1 \wedge (\theta-2) \wedge \epsilon$, we obtain
        \seqn{\label{11cpf}
        \frac{f_n(k/\sqrt{v\sigma^2 n};z_c)}{f_n(0;z_c)}
        & =
        \prod_{i=1}^n \left[1- v_i a\big(\frac{k}{\sqrt{v\sigma^2 n}}\big)
        + {\cal O}(\beta a\big(\frac{k}{\sqrt{v\sigma^2 n}}\big) i^{-\delta})
        \right]\\
        & =
        e^{-k^2/2d}[1+{\cal O}(k^{2+2\epsilon} n^{-\epsilon} )
        + {\cal O}(k^2 n^{- \delta})].}
With (\ref{hat{tau}_nlimit}), this gives the desired result.

\medskip \noindent {\em Proof of Theorem \ref{thm-1p}(b).}
Since $\delta<1 \wedge (\theta-2)$, it follows from (\ref{1.2b2}--\ref{1.2b3})
and (\ref{Delta_nlimit}--\ref{hat{tau}_nlimit}) that
        \eqn{
        \frac{\nabla^2f_n(0;z_c)}{f_n(0;z_c)}
        = -v\sigma^2  n [1+{\cal O}(\beta n^{-\delta})].
        }

\medskip \noindent {\em Proof of Theorem \ref{thm-1p}(c).} The claim is
immediate from Lemma~\ref{lem-Lpnorm}, which is now known to hold for all $n$.

\medskip \noindent {\em Proof of Theorem \ref{thm-1p}(d).}
Throughout this proof, we fix $z=z_c$ and drop $z_c$ from the notation.
The first identity of (\ref{eq:pthmd})
follows after we let $n \to \infty$ in (\ref{zetadef}), using
Lemma~\ref{zetan}.

To determine $A$, we use a summation argument. Let $\chi_n =
\sum_{k=0}^n f_k(0)$.  By (\ref{fkrec}),
    \seqn{
    \chi_n &= 1 + \sum_{j=1}^n f_j(0)
    = 1 + \sum_{j=1}^n \sum_{m=1}^{j}g_m(0) f_{j-m}(0) +  \sum_{j=1}^n e_j(0)
    \nonumber\\
    &= 1 + z\chi_{n-1}
    +\sum_{m=2}^{n}g_m(0)\chi_{n-m}+ \sum_{m=1}^ne_m(0).
    \label{S_nsum}}
Using (\ref{zetadef}) to rewrite $z$, this gives
    \eqn{
    f_n(0) = \chi_n - \chi_{n-1}
    = 1+ \zeta_n \chi_{n-1} - \sum_{m=2}^{n}
    g_m(0)(\chi_{n-1}-\chi_{n-m}) + \sum_{m=1}^ne_m(0).
    }
By Theorem~\ref{thm-1p}(a), $\chi_n \sim nA$ as $n \to \infty$.
Therefore, using Lemma~\ref{zetan} to bound the $\zeta_n$ term,
taking the limit $n \to \infty$ in the above equation gives
    \eqn{
    A = 1 - A \sum_{m=2}^\infty (m-1)g_m(0) + \sum_{m=1}^\infty e_m(0).
    }
With the first identity of (\ref{eq:pthmd}), this gives the second.

Finally, we use (\ref{Delta_nlimit}), (\ref{Delta_n}) and
Lemma~\ref{lem-pibds} to obtain
    \eqn{
    v = \lim_{n \to \infty}  v_n
    = \frac{-\sigma^{-2}\sum_{m=2}^\infty \nabla^2 g_m(0)}
    {1 + \sum_{m=2}^\infty (m-1) g_m(0)}.
    }
The result then follows, once we rewrite the denominator
using the first identity of (\ref{eq:pthmd}).

\subsection{Proof of Theorem~\protect\ref{thm-zc}}

By Theorem~\ref{thm-1p}(a), $\chi(z_c)=\infty$.  Therefore $z_c \geq z_c'$.
We need to rule out the possibility that $z_c>z_c'$.
Theorem~\ref{thm-1p} also
gives (\ref{fbdsp}) at $z=z_c$.  By assumption, the series
\eqn{
       G(z) = \sum_{m=2}^\infty g_m(0;z),
       \quad
       E(z) = \sum_{m=2}^\infty e_m(0;z)
}
therefore both converge
absolutely and are ${\cal O}(\beta)$ uniformly in $z \leq z_c$.
For $z<z_c'$, since the series defining $\chi(z)$ converges absolutely,
the basic
recursion relation (\ref{fkrec}) gives
\eqn{
       \chi(z) = 1 + z\chi(z) + G(z)\chi(z) +E(z),
}
and hence
\eqn{
\label{chiEG}
       \chi(z) = \frac{1+E(z)}{1-z-G(z)}, \quad (z<z_c').
}
It is implicit in the bound on $\partial_z g_m(k;z)$ of Assumption~G
that $g_m(k;\cdot)$ is continuous on $[0,z_c]$.
By dominated convergence, $G$ is
also continuous on $[0,z_c]$.
Since $E(z) = {\cal O}(\beta)$ and
$\lim_{z \uparrow z_c'}\chi(z) = \infty$, it then follows from (\ref{chiEG}) that
\eqn{
\label{pc0}
       1 - z_c' - G(z_c') = 0.
}
By the first identity of (\ref{eq:pthmd}), (\ref{pc0}) holds also when $z_c'$ is replaced by $z_c$.
If $z_c' \neq z_c$, then it follows from the mean-value theorem that
\eqn{
       z_c - z_c' = G(z_c') - G(z_c) = -(z_c-z_c') \sum_{m=2}^\infty
       \partial_z g_m(0;t)
}
for some $t \in (z_c',z_c)$.  However, by a bound of Assumption~G, the sum
on the right side is ${\cal O}(\beta)$ uniformly in $t \leq z_c$.
This is a contradiction, so we conclude that $z_c=z_c'$.
\qed

\section*{Acknowledgements}
A version of this work appeared in the PhD thesis \cite{Holm05}.  The work of RvdH and MH was supported in part by Netherlands Organisation for Scientific Research (NWO).  The work of GS was supported in part by NSERC of Canada.


\begin{thebibliography}{1}

\bibitem{HHS07}
R.~van~der Hofstad, M.~Holmes, and G.~Slade.
\newblock An extension of the generalised inductive approach to the lace
  expansion.
\newblock Preprint, 2007.

\bibitem{HS02}
R.~van~der Hofstad and G.~Slade.
\newblock A generalised inductive approach to the lace expansion.
\newblock {\em Probab. Theory Relat. Fields.}, 122:389--430, 2002.

\bibitem{H07}
M.~Holmes.
\newblock Convergence of lattice trees to super-{B}rownian motion above the
  critical dimension.
\newblock Preprint, 2007.

\bibitem{Holm05}
M.~Holmes.
\newblock {\em Convergence of lattice trees to super-{Brownian} motion above
  the critical dimension}.
\newblock PhD thesis, University of British Columbia, (2005).


\end{thebibliography}

\end{document}